# ABOUT STABLE PERIODIC HELIXES, L-ITERATION AND CHAOS GENERATED BY UNBOUNDED FUNCTIONS

## UNPREDICTABILITY, APERIODICITY, SENSITIVE DEPENDENCE AND CHAOTIC BEHAVIOR

by Andrei Vieru

**Abstract**

We consider *stable periodic helixes* as a generalization of stable periodic orbits. We see that in the studied class of iterated functions Chaos always arise *suddenly*. Therefore, we'll study the route from chaos to order rather than the route from order to chaos. We show that, paradoxically, genuine Chaos may look as much like Order and during as many iteration steps as one may wish. The way toward order through the lengthening of the order-like iteration steps series brings to the fore a constant that permits to calculate precisely the *order point*, i.e. the parameter value for which chaos turns into order.

Then, we'll propose a generalization of the idea of map iteration that do not imply the existence of periodic orbits. We'll show that, within a strictly deterministic context, unpredictability, 'aperiodic order', sensitive dependence and chaos are completely different concepts and we'll try to show what this difference is made of. We'll also propose an example of non-chaotic 'aperiodic order'.

## 1. Unbounded functions, periodic helixes and chaos

It was often said that Chaos in general implies functions we deal with are bounded and the metric spaces compact. We consider this idea as a foregone conclusion, whose origins are probably deeply rooted in mathematicians' attachment to the strict concept of periodic (or quasi-periodic) 'orbit'.

We'll show later that the concept of periodic helix can be applied to bounded and unbounded functions as well.

**Definition 1**

Let $\{u(n)\}_{n \in N}$ be a sequence. Assume there is an integer $j$ and a real $r$ such as
$\forall i \in N \; \forall k \in N \; u(i + kj) - u(i) = mr$ (where $m$ is an integer[1])
We'll call the sequence $\{u(n)\}_{n \in N}$ a **helix** with period $j$ and with modulo $r$.
Examples:
a) 1, 2, 1, 2, 1, 2, … is a helix with period 2 and modulo 0. (So, obviously any sequence based on periodicity – i.e. on an orbit – is a helix modulo 0.)
b) the sequence 1, 1, 2, 3, 5, 8, 13, 21, 34,… of the Fibonacci numbers is a helix with period 3 and modulo 2. It is also a helix with period 8 and modulo 3. It is also a helix with period 6 and modulo 4, etc.

---

[1] one can distinguish the case where $m$ is a constant from the case when it is not; one can thus distinguish helixes from spirals or other varieties. We'll not examine here these details.

c) Trivially, *N* is a helix with period 1 and modulo 1. It is also a helix with period 2 and modulo 2, etc.

Let us consider the following sequence: $u(1)=0.5$ $\forall n \geq 1$ $u(n+1)=0.4\sin(\pi u(n))+u(n)+1$. It is easy to see that this sequence is 'attracted' by *N*, because $\forall \varepsilon > 0$ $\exists N(\varepsilon)$ $\forall n \geq N(\varepsilon)$ $|u(n-1) - n| < \varepsilon$. We'll say that this sequence is attracted by a stable helix with period 1 (we'll consider, from now on, only the smallest period of a helix).

Consider the family of sequences $u_a(1)=a$ $\forall n$ $u_a(n+1)=0.4\sin(\pi u(n))+u(n)+b$ where $0.88725599 < b < 1$. Calculations reveal that the sequences are attracted by a stable periodic helix with period 2 and with modulo 2. (So 1 was a bifurcation point.)

Let's see what happens when we iterate functions of the form $f_b(x)= 0.4\sin(\pi x)+x+b$ and when the *b* value decreases.

Near $b = 0.88725599...$, we are about to sink into chaos. The slide into chaos does not follow the Feigenbaum scenario. We'll try to briefly describe this new way of *sudden* skid into chaos bellow and, more detailed, in APPENDIX 3.

Something interesting seems to happen iterating functions, setting the parameter *b* near 0.88725598... Let's consider the family of sequences:

$u_a(1)=a$ $\forall n$ $u_a(n+1) = 0.4\sin(\pi u_a(n)) + u_a(n) + 0.88725598$ **(1)**

The behavior of the sequences seems 'almost periodic'. In fact, strictly speaking, it isn't periodic. Actually, we find instead of a helix, a *pseudo-helix*: the 'would-be helix' (which, considering an insufficiently long part of an infinite sequence, may seem to be a period 2 periodic helix[2]) undergoes some 'slithering'... From time to time – quasi periodically[3] – it escapes the 'would-be (period 2) helix', takes '*chaotic*' values faraway from it, then comes close again to the 'would-be helix', several hundreds – or thousands – of terms later escapes again, takes again faraway, *chaotic* – but different! - values from the 'would-be helix', and so on...

Following our[4] definition of chaos, the family of sequences **(1)** is already chaotic, since it should be not so difficult to show that it satisfies the chaos criteria[5], in particular *sensitive dependence on initial conditions*.

Consider the family of sequences $u_a(1) = a$ $\forall n$ $u_a(n+1) = 0.4\sin(\pi u_a(n)) + u_a(n) + 0.882$. Calculations reveal chaotic behavior.

---

[2] We'll not give here any definition of quasi-periodic helix, since the analogy with the concept of quasi-periodic orbit make it obvious.

[3] the *average order* of the quasi periodicity depends on *a*.

[4] We daresay 'our' because, unlike all other authors, we suppose neither compactness of the metric space nor boundedness of the functions.

[5] to say it more precisely: 'We believe that chaos should involve not only nearby points can diverge apart but also faraway points can get close to each other. Therefore, we propose to call a continuous map *f* from an infinite [not necessarily] compact metric space (X, d) to itself chaotic if there exists a positive number λ such that for any point *x* and any nonempty open set V (not necessarily an open neighborhood of *x*) in X there is a point *y* in V such that $\lim \sup_{n \to \infty} d(f^n(x), f^n(y)) \geq \lambda$ and $\lim \inf_{n \to \infty} d(f^n(x), f^n(y)) = 0$.' (Bau-Sen Du *'On the nature of chaos'* arXiv:math.DS/0602585 v1)

(We added the red words in the brackets; they miss in Bau-Sen Du's definition, whose point of view doesn't entirely coincides with ours. See also our article *'General definitions of chaos for continuous and discrete-time processes'* arXiv:0802.0677 [math.DS] §4.2

For $b = 0.8811$, $b = 0.881$ and $b = 0.88093$ we find that all sequences are attracted by stable periodic helixes of order (period) 9, 18 and 36 respectively.

For $b = 0.87$ and $b = 0.8698$ we find that all sequences are attracted by periodic helixes of order 7 and 14 respectively.

For $b = 0.801$, $b = 0.799$ and $b = 0.7945$ we find that all sequences are attracted by periodic helixes of order 20, 10 and 5 respectively[6]. (See APPENDIX 1.)

For $b \in [0.7753369\ldots, 0.6635\ldots]$ we find stable periodic helixes of order 3, which attract the sequences. (See APPENDIX 2.)

Needless to say, *sensitive dependence* and *metric independence with respect to initial conditions* disappear as soon as stable periodic helixes appear and vice-versa.

The initial stable periodic helixes of order 1, the stable helixes of order 2 and stable helixes with period 9, 7, 5 and 3 beyond chaos point seem to be a part of the Sharkovsky's ordering, as it happens with periodic orbits that appear before and beyond chaos point reached in bifurcation processes[7]. But only a part: there are no periodic helixes of order $2^n$ with $n>1$. Detecting the existence – or the proof of the inexistence – of periods of order 11, 13, 15, 17 and so on are tasks that we propose to the reader;

**Definition 2**

*Topological transitivity modulo* 1. Topological transitivity modulo 1 means that given an arbitrary small subinterval S of [0, 1] we'll find, for every initial $x$, an iterate whose fractional part belongs to S. (The general definition of *Topological transitivity modulo r* is proposed to the reader.)

We also propose to the reader to prove that for 'chaotic' $b$ values the family of sequences we are speaking about has topological transitivity modulo 1.

So far we haven't studied the question of the existence of *instable* (i.e. non attracting) periodic helixes and even less the question of their density, if they exist. Again we suggest the reader to examine it.

Without loss of generality, we'll formulate an additional question in relationship with a concrete numerical example:

The family of sequences generated iterating $f_{0.41}(x) = 0.4\sin(\pi x)+x+0.41$ has, conjecturally, chaotic behavior (in the sense of footnote 5 from which the word 'compact' should be either withdrawn or preceded by the 'not necessarily'). To illustrate it, we give the following result:

$u_{0.5}(1,000,000) = 99081.47278\ldots$ and $u_{0.500001}(1,000,000) = 99025.98068\ldots$
but
$u_{0.5}(100,000) = 9893.52525\ldots$ and $u_{0.500001}(100,000) = 9903.48937\ldots$

---

[6] Curiously, as $a$ decreases from 0.801 to 0.7945, we don't have period-doubling, but 'period-dividing' by 2.

[7] We only notice that the interval in which the parameter $a$ value corresponds to a helix of order 3 is the longest one, exactly as it happens to be the case in classical bifurcation processes (for example, in the logistic map).

1°) Does the sign of $u_{0.5}(n)-u_{0.500001}(n)$ change infinitely many times[8]?

2°) Is the sequence $|u_{0.5}(n)-u_{0.500001}(n)|$ bounded or not? More generally, is the sequence $|u_\alpha(n)-u_{\alpha+\varepsilon}(n)|$ bounded? If yes, how does its upper bound depend on $\alpha$ and on $\varepsilon$?

3°) What is the inferior bound? (We state that the inferior bound is 0: $\forall a \forall \varepsilon \forall \delta \exists n \; |u_a(n) - u_{a+\delta}(n)| < \varepsilon$. Can the reader prove it?)

## 2. Unpredictability, aperiodic order, sensitive dependence and chaos

'Aperiodic order' does not imply chaos, because it does not imply sensitive dependence on initial conditions.
Example:
Consider one of the simplest Lindenmayer systems (we'll designate it by $L_1$), based on the alphabet {A, B} and on the rules A→AB and B→BA. Beginning from the one-letter word A, 'simultaneous' endless iteration of these rules produces the infinite word ABBABAABBAABABBA…, whose features, including aperiodicity, are well known. We'll designate this word by $<A, B>^\infty$. Assume A and B represent functions.

Let's take as an example $A(x) = \Gamma(x+1)$ and $B(x) = \cos x$.

Consider the family of sequences $(U_a(n))_{a \in R}$ constructed in the following way:
$U_a(0) = a$
$U_a(1) = A(a) = \Gamma(a+1)$
$U_a(2) = (B \circ A)(a) = \cos(\Gamma(x+1))$
$U_a(3) = (B \circ B \circ A)(a) = \cos(\cos(\Gamma(a+1)))$
$U_a(4) = (A \circ B \circ B \circ A)(a)$, $U_a(5) = (B \circ A \circ B \circ B \circ A)(a)$
 etc.
where in the $n$-th term of the sequence we have, reading from right to left, the first $n$ letters (from left to right) of the infinite word $<A, B>^\infty$.

For any given $a$, the sequence $U_a(n)$ is aperiodic and, *in some sense*, unpredictable: there is no general formula of its $n$-th term, as there is no general formula for the $n$-th letter of the word $<A, B>^\infty$.

Any of these sequences – and even the entire family itself – may be considered as an example of 'aperiodic order'. Their aperiodicity is interesting, there are features that in my opinion deserve to be studied[9], but this family cannot be considered as an example of 'chaos' or 'chaotic behavior'. For a very simple reason: there is no Sensitive Dependence on Initial Conditions (SDIC) here! There is neither *density of orbits* (there are no orbits at all) nor *topological transitivity*!

If we know one of the sequences of this family, we can 'predict' the behavior of all of them: they attract each other. To say it in a precise way:

$\forall a \; \forall b \; \forall \varepsilon \; \exists N(a, b, \varepsilon) \in N \; \forall n \geq N(a, b, \varepsilon) \; |u_a(n) - u_b(n)| < \varepsilon$

---

[8] Or, more generally, if the sign of $u_a(n)-u_{a+\varepsilon}(n)$ changes infinitely many times
[9] For example, there are no orbits, but there is something we would like to call 'water slides', a concept we shall not define now. (The reader is proposed to do it.)

Let's now consider a family of sequences $(v_a(n))$ ($a \in \mathbf{R}$) constructed in a similar way, only we'll write [10]  $A(x)=\Gamma(x)$ if $x \notin \mathbf{Z}\text{-}$ and if $x \neq 0$
$A(x)=\Gamma(x+1/2)$ if $x \in \mathbf{Z}\text{-}$ or if $x=0$
and $B(x)=\sin(\pi x)$

It is easy to see that the family of sequences $(v_a(n))$ has the SDIC property. It is also easy to see that this family contains unbounded sequences. It is very likely it is chaotic (in the sense of the footnote 5 from which the word 'compact' should be withdrawn). If metric independence were proved, then, according to this modified definition[11], the family of sequences $v_a(n)$ might be considered as 'chaotic'.

## 3. L-iteration, unbounded functions and chaos

### Definition 3
In the context of unbounded functions topological transitivity may be defined as usually: a function $f$ is topologically transitive if, given any two intervals **U** and **V**, there is some positive integer $k$ such that $f^k(\mathbf{U}) \cap \mathbf{V} \neq \emptyset$.

A finite set of functions $\{f_1,\ldots, f_n\}$ is topologically transitive with respect to a Lindenmayer system **L** based on the alphabet $\{f_1,\ldots, f_n\}$ if, given any two intervals **U** and **V**, there is some positive integer $k$ such that $f_\mathbf{L}^k(\mathbf{U}) \cap \mathbf{V} \neq \emptyset$, where $f^k_\mathbf{L}(\mathbf{U})= \{f_\mathbf{L}^k(x) | x \in \mathbf{U}\}$, where $f_\mathbf{L}^k(x)$ designates $(f_{\mathbf{L}(k)} \circ \ldots \circ f_{\mathbf{L}(2)} \circ f_{\mathbf{L}(1)})(x)$, where $f_{\mathbf{L}(1)}, f_{\mathbf{L}(2)},\ldots, f_{\mathbf{L}(k)}$ are, in their exact order, the first $k$ letters of the infinite word generated by **L**.

More generally, a family of sequences $w_a(n)$ (where $a \in \mathbf{D} \subset \mathbf{R}$ and such as $w_a(1)=a$) is topologically transitive if, given any two intervals $\mathbf{U} \subset \mathbf{D}$ and $\mathbf{V} \subset \mathbf{D}$, there is some positive integer $k$ and some $\lambda \in \mathbf{U}$ such that $w_\lambda(k) \in \mathbf{V}$.

Example 1: { $A(x)=\Gamma(x)$ if $x \notin \mathbf{Z}\text{-}$ and if $x \neq 0$  $A(x)=\Gamma(x+1/2)$ if $x \in \mathbf{Z}\text{-}$ or if $x=0$, $B(x)=\sin\pi x$} is (likely to be) topologically transitive with respect to $L_1$

Example 2: $\{A(x)=\Gamma(x+1), B(x)=\cos x\}$ is definitely not topologically transitive with respect to $L_1$.

### Definition 4
A family of functions $\{f_a\}_{a \in \mathbf{R}+}$ has the helix density property if there is a (somewhere) dense collection of points engendering helixes. We suggest the reader to examine this property and the objects it is related to.

---

[10] the fact that $A(x)$ is defined as discontinuous function shouldn't bother anyone. In fact we'll consider a family of sequences, which are, of course, continuous functions. It is easy to construct chaotic families of sequences based on L-iteration of continuous functions depending on some parameters: for example $A(x)= a\cos x+b$ and $B(x)=cx^2$
We confess we have no idea if for some values of the parameters these families depend on they might have a non-chaotic behavior.

[11] We have, by the way, to ask the following question. How should we call a family of sequences $(u_a(n))$ that satisfies the following condition: there exists a positive number $\lambda$ such that for any point $\alpha$ and any nonempty open set V (not necessarily an open neighborhood of $\alpha$) in $\Omega$ there is a point $\beta$ in V such that $\lim \sup_{n \to \infty} d(u_\alpha(n), u_\beta(n)) = \infty$ and $\lim \inf_{n \to \infty} d(f^n(x), f^n(y)) \leq \lambda$? (We formulate this question using words and notations borrowed from Bau-Sen Du's article *On the Nature of Chaos*, available on arXiv)

*Conclusion*

The extension of the concept of periodic orbit to the concept of periodic helix, the extension of the concept of iteration to L-iteration, along with our chaos definitions[12] permit to see that families of unbounded real functions can generate chaos.

We have chosen examples – see also the appendixes – in which the fractional part of the terms of the sequence is periodic, while the integer part displays monotone growth. We attribute some philosophical importance to the idea that (quasi) periodicity and/or pseudo (quasi) periodicity may be reached on some levels, while on some distinct levels there is not even pseudo periodicity at all. We consider that even generalizations without practical relevance are important if one wants to formulate correct general definitions.

In our text, a lot of questions remain unanswered: we propose the reader to try to give an answer. An important question has already been raised at the end of our paper *'General definitions of chaos for continuous and discrete-time processes'*. We repeat it here, because of the unfathomable variety of chaotic behaviors. Are these definitions recursive? Once again, we are rather pessimistic: we believe there cannot exist any general recursive definition of chaos.

**Acknowledgements**

We express our deep gratitude Vlad Vieru and Dmitry Zotov for computer programming and to Robert Vinograd for useful discussions.

Andrei Vieru
andreivieru2004@yahoo.fr

**References**

(1) Henri Poincaré *Sur le problème des trois corps et les équations de la dynamique*, Acta Mathematica **13** (1890), 1–270.
(2) Bau-Sen Du *On the nature of chaos*, arXiv:math/0602585 (2006)
(3) Vladimir Arnold *Geometry and Dynamics of Galois Fields*, UMN 2004
(4) Vladimir Arnol'd *Ordinary Differential Equations*
(5) Sharkovsky, A.N. [1964] : *Coexistence of cycles of continuous mapping of the line into itself*. Ukrainian Math. J., **16**, 61-71 [in Russian]; [1995] : Intern. J. Bifurcation and Chaos **5**, no. 5, 1263-1273.
(6) Andrei Vieru *Generalized Iteration, Catastrophes, Generalized Sharkovsky's ordering* arXiv:0801.3755 [math.DS]
(7) Andrei Vieru *General definitions of chaos for continuous and discrete-time processes* arXiv:0802.0677 [math.DS]
(8) Andrei Vieru, *Le gai Ecclésiaste,* Editions du Seuil, 2007, p.194-195

---

[12] See our paper *'General definitions of chaos for continuous and discrete-time processes'* arXiv:0802.0677 [math.DS]

**APPENDIX 1 (an example of a stable periodic helix of order 10)**

$u_a(1)=0.5$ $\forall n$ $u_a(n+1)=0.4\sin(\pi u_a(n))+u_a(n)+0.8$

Term number 1000000: 800000.93556582…
Term number 1000001: 800001.81598435…
Term number 1000002: 800002.39741014…
Term number 1000003: 800003.57681447…
Term number 1000004: 800003.98840510…
Term number 1000005: 800004.77383775…
Term number 1000006: 800005.83472582…
Term number 1000007: 800006.43624316…
Term number 1000008: 800007.62824609…
Term number 1000009: 800008.06027449…

Term number 1000010: 800008.93556582…
Term number 1000011: 800009.81598435…
Term number 1000012: 800010.39741014…
Term number 1000013: 800011.57681447…
Term number 1000014: 800011.98840510…
Term number 1000015: 800012.77383776…
Term number 1000016: 800013.83472582…
Term number 1000017: 800014.43624316…
Term number 1000018: 800015.62824609…
Term number 1000019: 800016.06027449…

Term number 1000020: 800016.93556582…
Term number 1000021: 800017.81598435…
Term number 1000022: 800018.39741014…
Term number 1000023: 800019.57681446…
Etc…

**APPENDIX 2 (an example of a stable periodic helix of order 3)**

$u_a(1)=0.8$ $\forall n$ $u_a(n+1)=0.4\sin(\pi u_a(n))+u_a(n)+0.7$

| Term nr | | |
|---|---|---|
| | 20000 | 13 333.7162148952… |
| | 20001 | 13 334.1049995507… |
| | 20002 | 13 334.9345659839… |
| | 20003 | 13 335.7162148952… |
| | 20004 | 13 336.1049995507… |
| | 20005 | 13 336.9345659839… |
| | 20006 | 13 337.7162148952… |
| | 20007 | 13 338.1049995507… |
| | 20008 | 13 338.9345659839… |
| | 20009 | 13 339.7162148952… |
| | 20010 | 13 340.1049995507… |
| | 20011 | 13 340.9345659839… |
| | 20012 | 13 341.7162148952… |
| | 20013 | 13 342.1049995507… |
| | 20014 | 13 342.9345659839… |
| | 20015 | 13 343.7162148952… |
| | 20016 | 13 344.1049995507… |
| | 20017 | 13 344.9345659839… |

**APPENDIX 3: THE ROUTE FROM CHAOS TO ORDER**

Consider the family of sequences
$u_a(1)=a \quad \forall n \quad u_a(n+1)=0.4\sin(\pi u_a(n))+u_a(n)+0.8872559$

*pseudo*[13]*-quasi-periodic helixes* of order 2 seam to appear. Here is, for *a*=0.5, a sample easy to interpret if one starts to pay attention to the first seven or eight digits in the fractional part and then looks to the rest. The pseudo-helix slowly slides down:

term nr. :
| | |
|---|---|
| 1200 | 1 198.6422740820262334… |
| 1201 | 1 199.8902348292253919… |
| 1202 | 1 200.6422732472858570… |
| 1203 | 1 201.8902344478810846… |
| 1204 | 1 202.6422724149413170… |
| 1205 | 1 203.8902340676286258… |
| 1206 | 1 204.6422715849801079… |
| 1207 | 1 205.8902336884627857… |
| 1208 | 1 206.6422707573906337… |
| 1209 | 1 207.8902333103778801… |
| 1210 | 1 208.6422699321606160… |
| 1211 | 1 209.8902329333684520… |
| 1212 | 1 210.6422691092784589… |
| 1213 | 1 211.8902325574292718… |
| 1214 | 1 212.6422682887321116… |
| 1215 | 1 213.8902321825548825… |
| 1216 | 1 214.6422674705104328… |
| 1217 | 1 215.8902318087400545… |
| 1218 | 1 216.6422666546015989… |
| 1219 | 1 217.8902314359795582… |
| 1220 | 1 218.6422658409937867… |

etc.

Nevertheless, every approximately 6090 terms, the sequence completely escapes any kind of 'approximate order 2 periodicity' in the fractional part: it turns out that, *in the strict sense*, there is no helix here. Here are the first four such *obviously chaotic episodes*[14] *(skids)* in this sequence:

term nr:
| | |
|---|---|
| 6090 | 6 088.5355653501092092… |
| 6091 | 6 089.8203270456033351… |
| 6092 | 6 090.4935993405651971… |
| 6093 | 6 091.7807743748289795… |
| 6094 | 6 092.4138112253767758… |
| 6095 | 6 093.6864932140988458… |
| 6096 | 6 094.2404600462305098… |
| 6097 | 6 095.4019559273619961… |
| 6098 | 6 095.9080368775057650… |
| 6099 | 6 096.6813294927551397… |
| 6100 | 6 097.9054184193000765… |
| 6101 | 6 098.6755608477915303… |
| 6102 | 6 099.9035040366752582… |
| 6103 | 6 100.6713483331050156… |

etc.

One should note that at the end – the last four or six terms – show that the sequence is 'picked' up by what might again seem as a *pseudo-quasi-periodic* helix. In fact it looks again almost like a helix during some other 6090 terms; then, the sequence 'escapes from' the *pseudo-quasi-periodic helix* again:

term nr:
| | |
|---|---|
| 12170 | 12 167.8574588992596546… |
| 12171 | 12 168.5715192090974597… |
| 12172 | 12 169.8487209133254510… |
| 12173 | 12 170.5529499261192541… |
| 12174 | 12 171.8346843046238064… |
| 12175 | 12 172.5234122360525362… |

---

[13] For example, for *b* = 0.88725598 the sequence escapes for the first time from the pseudo-quasi-periodic helix after the 18900th term. The *average order of periodicity of the obviously chaotic episodes (skids)* does not depend on the initial *x* value.

[14] Of course it is correct to consider the whole family of sequences as chaotic and not separately one concrete sequence or, even less, only the *obviously chaotic episodes (skids)*.

| | | |
|---|---|---|
| 12176 | 12 173.8095866529656632… | |
| 12177 | 12 174.4715797833596298… | |
| 12178 | 12 175.7572423891324433… | |
| 12179 | 12 176.3681636470391823… | |
| 12180 | 12 177.6215988233689131… | |
| 12181 | 12 178.1376884313867777… | |
| 12182 | 12 179.1926232361438451… | |
| 12183 | 12 179.8523270162313565… | |
| 12184 | 12 180.5605970402139064… | |
| 12185 | 12 181.8406265641915525… | |
| 12186 | 12 182.5358713395726227… | |
| 12187 | 12 183.8205899778731691… | |
| 12188 | 12 184.4941415023167792… | |
| 12189 | 12 185.7813296553249529… | |
| 12190 | 12 186.4149056262740487… | |
| 12191 | 12 187.6879531830654741… | |
| 12192 | 12 188.2429379649311159… | |
| 12193 | 12 189.4066923255595611… | |
| 12194 | 12 189.9110111082736694… | |
| 12195 | 12 190.6789912722550012… | |
| 12196 | 12 191.9074615991958126… | |
| 12197 | 12 192.6800614462918020… | |
| 12198 | 12 193.9050071359652065… | |
| 12199 | 12 194.6746554767996713… | |
| 12200 | 12 195.9031934535796609… | |
| 12201 | 12 196.6706653136970999… | |
| 12202 | 12 197.9017916846714797… | |
| 12203 | 12 198.6675840349471400… | |
| 12204 | 12 199.9006722591439029… | |
| 12205 | 12 200.6651250754603097… | |
| 12206 | 12 201.8997556809663365… | |
| 12207 | 12 202.6631128259468824… | |
| 12208 | 12 203.8989901950917556… | |
| 12209 | 12 204.6614330703250744… | |
| 12210 | 12 205.8983405271283118… | |
| 12211 | 12 206.6600080275693472… | |
| 12212 | 12 207.8977817396680621… | |
| 12213 | 12 208.6587827490275231… | |
| 12214 | 12 209.8972956672496366… | |
| 12215 | 12 210.6577172324032290… | |
| 12216 | 12 211.8968687423675874… | |
| 12217 | 12 212.6567816166825651… | |
| 12218 | 12 213.8964906143519329… | |
| 12219 | 12 214.6559531316906941… | |
| 12220 | 12 215.8961532412304223… | |

etc.

One should note that that sequence is again 'picked up' by the *pseudo-helix*. It will escape from it again about 6090 terms later:

| | | |
|---|---|---|
| 18267 | 18 262.5820201598580752… | |
| 18268 | 18 263.8560701999231242… | |
| 18269 | 18 264.5685591392903007… | |
| 18270 | 18 265.8465727220536792… | |
| 18271 | 18 266.5484055959532270… | |
| 18272 | 18 267.8310453048143245… | |
| 18273 | 18 268.5158163782434713… | |
| 18274 | 18 269.8025785880854528… | |
| 18275 | 18 270.4573495707300026… | |
| 18276 | 18 271.7410201608145144… | |
| 18277 | 18 272.3375676760188071… | |
| 18278 | 18 273.5738635253583197… | |
| 18279 | 18 274.0718405451661965… | |
| 18280 | 18 275.0486094656698697… | |
| 18281 | 18 275.8750180548522621 | (here there is an 'attempt' to hang on the pseudo- |
| 18282 | 18 276.6092215434982791 | helix. But the sequence doesn't hang high |
| 18283 | 18 277.8731599913480750 | enough, so it quickly slides down again a |
| 18284 | 18 278.6052088721880864 | dozen of terms later) |

| | |
|---|---|
| 18285 | 18 279.8708138135298213… |
| 18286 | 18 280.6001496339595178… |
| 18287 | 18 281.8677699922591273… |
| 18288 | 18 282.5935988221935986… |
| 18289 | 18 283.8636859603320772… |
| 18290 | 18 284.5848325500046485… |
| 18291 | 18 285.8579668891870824… |
| 18292 | 18 286.5726028356111783… |
| 18293 | 18 287.8494988904349157… |
| 18294 | 18 288.5545977375440998… |
| 18295 | 18 289.8359839631302748… |
| 18296 | 18 290.5261313903611153… |
| 18297 | 18 291.8120401561864128… |
| 18298 | 18 292.4765877246900345… |
| 18299 | 18 293.7627621379797347… |
| 18300 | 18 294.3787397266787593… |
| 18301 | 18 295.6373203064431436… |
| 18302 | 18 296.1612244574207580… |
| 18303 | 18 297.2425287698315515… |
| 18304 | 18 297.8536580154795956… |
| 18305 | 18 298.5634253934877052… |
| 18306 | 18 299.8427668807125883… |
| 18307 | 18 300.5403754402505001… |
| 18308 | 18 301.8244178135173570… |
| 18309 | 18 302.5020508484740276… |
| 18310 | 18 303.7892984462341701… |
| 18311 | 18 304.4306955213432957… |
| 18312 | 18 305.7085078551790502… |
| 18313 | 18 306.2785559746007493… |
| 18314 | 18 307.4728573245592997… |
| 18315 | 18 307.9615665801320574… |
| 18316 | 18 308.8006428854532714… |
| 18317 | 18 309.9223588237437070… |
| 18318 | 18 310.7130125339499500… |
| 18319 | 18 311.9139960599459300… |
| 18320 | 18 312.6944863926182734… |
| 18321 | 18 313.9093728721963998… |
| 18322 | 18 314.6842757662525401… |
| 18323 | 18 315.9063534101151163… |
| 18324 | 18 316.6776197987601336… |

and again:

| term nr. | |
|---|---|
| 24398 | 24 390.5378065753138799… |
| 24399 | 24 391.8222443921222293… |
| 24400 | 24 392.4975562131548941… |
| 24401 | 24 393.7848003247709130… |
| 24402 | 24 394.4217633864100208… |
| 24403 | 24 395.6969976839754963… |
| 24404 | 24 396.2584436047327472… |
| 24405 | 24 397.4359446217931691… |
| 24406 | 24 397.9312724046503718… |
| 24407 | 24 398.7328321440036234… |
| 24408 | 24 399.9177670327626402… |
| 24409 | 24 400.7028315750321781… |

etc.

        One should also note that the length of these obviously chaotic episodes is not always the same. As we already said, it is also easy[15] to notice – choosing close *x* values – that sensitive dependence on initial conditions appears as soon as stable periodic helixes disappear.

---

[15] In fact, much more difficult if the parameter *b* value (the free term) is very close to the *order point* (where pseudo helixes become true helixes and chaos becomes order). As the *average order of the periodicity of the obviously chaotic episodes* decreases under about 90, the sensitive dependence is easier to detect experimentally, but the periodicity itself much harder to notice, because the border between what is *obviously chaotic* and what looks more or less like a *pseudo-quasi-periodic helix* tends to vanish.

Choosing lower *b* values we may find lower (quasi) periodicities of the subsequences in which we find terms that escape from the *pseudo-quasi-periodic helixes*. Choosing higher *b* values we may find higher (probably as high as we wish) (quasi) periodicities of the *obviously chaotic episodes* (skids).

An example of lower pseudo-quasi-periodicity of the *obviously chaotic episodes*:

$$u_a(1)=0.8 \quad \forall n \quad u_a(n+1)=0.4\sin(\pi u_a(n))+u_a(n)+0.8872$$

One can see that the subsequences made of red numbers, when the sequence escapes from the *pseudo-quasi-periodic helix* of order 2 lead again to a *pseudo-quasi-periodic helix* of order 2. The subsequences made of red numbers – the *obviously chaotic episodes* – appear approximately every 260 or 270 terms of the sequence (for *b*=0.8872)

Term nr:
| | |
|---|---|
| 19626 | 19 479.8703111548493325… |
| 19627 | 19 480.5990109244921769… |
| 19628 | 19 481.8670157738633861… |
| 19629 | 19 482.5919219877469004… |
| 19630 | 19 483.8625586331581871… |
| 19631 | 19 484.5823616535199108… |
| 19632 | 19 485.8562462140507705… |
| 19633 | 19 486.5688782371835259… |
| 19634 | 19 487.8467500226033735… |
| 19635 | 19 488.5487244431024010… |
| 19636 | 19 489.8312473572477757… |
| 19637 | 19 490.5161815427818510… |
| 19638 | 19 491.8028647980463575… |
| 19639 | 19 492.4578726503314101… |
| 19640 | 19 493.7415746164770098… |
| 19641 | 19 494.3385452556976816… |
| 19642 | 19 495.5753835831710603… |
| 19643 | 19 496.0737484238561592… |
| 19644 | 19 497.0527965305045655… |
| 19645 | 19 497.8739542460898519… |
| 19646 | 19 498.6068676007271279… |
| 19647 | 19 499.8717350343467842… |
| 19648 | 19 500.6020792162817088… |
| 19649 | 19 501.8688863068382489… |
| 19650 | 19 502.5959437204983260… |
| 19651 | 19 503.8651105400713277… |
| 19652 | 19 504.5878314084438898… |
| 19653 | 19 505.8599002502633084… |
| 19654 | 19 506.5766751226401539… |
| 19655 | 19 507.8523262997914571… |
| 19656 | 19 508.5605396186329017… |
| 19657 | 19 509.8405268902206444… |
| 19658 | 19 510.5356058955330809… |
| 19659 | 19 511.8203060068044579… |
| 19660 | 19 512.4935000652367307… |
| 19661 | 19 513.7806166716509324… |
| 19662 | 19 514.4134446495772863… |
| 19663 | 19 515.6859472703508800… |
| 19664 | 19 516.2394793483072135… |
| 19665 | 19 517.4000208821416891… |
| 19666 | 19 517.9067901674461609… |
| 19667 | 19 518.6785260291944724… |
| 19668 | 19 519.9044460581571911… |
| 19669 | 19 520.6733647753026162… |
| 19670 | 19 521.9026900098615442… |
| 19671 | 19 522.6695025068875111… |
| 19672 | 19 523.9013171335973311… |
| 19673 | 19 524.6664855392155005… |
| etc. | |